\documentclass[11pt]{article}
\usepackage[a4paper,margin=2.4cm]{geometry}
\usepackage{amsmath,amssymb,booktabs,array,enumitem}
\usepackage[english]{babel}

\newtheorem{lemma}{Lemma}

\newtheorem{theorem}[lemma]{Theorem}

\newcommand{\eps}{\varepsilon}
\newcommand{\e}{\varepsilon}
\newcommand{\Rd}{\mathbb{R}^d}
\newcommand{\sig}{\sigma_{d-1}}
\newcommand{\ac}{\overline a^{\rm c}}
\title{Homogenization of Gagliardo seminorms: the global picture}
\author{Andrea Braides, Stefano Mannella, Marco Picerni,\\ and Alec Jacopo Almo Schiavoni Piazza}
\date{}
\begin{document}
\maketitle

\begin{abstract}
We describe the asymptotic behaviour of squared Gagliardo seminorms with
periodically oscillating coefficients when the period of the oscillations
vanishes and, at the same time, the fractional exponent tends to a limit in $[0,1]$, the most relevant cases being zero or one. For a fixed exponent no homogenization effect
occurs and the coefficient is simply replaced by its average, while at the
endpoints the limit depends on the mutual rate of convergence of the two
parameters. We review and compare two recent results: a homogenized version
of the Maz'ya--Shaposhnikova formula, where the limit is a multiple of the
squared Lebesgue norm obtained by optimizing oscillations in the interactions
at finite and at infinite scale, and a homogenized version of the
Bourgain--Brezis--Mironescu formula, where the limit is a Dirichlet-type
energy whose matrix is a convex combination of the average of the coefficient
and of the homogenized matrix of its diagonal values. We then propose a
unified statement in terms of an expansion by Gamma-convergence, which
approximates the oscillating energies by homogeneous ones in all regimes.


\end{abstract}

\section{Introduction}
In this paper, we provide a description of the asymptotic behaviour of Gagliardo seminorms in the presence of oscillations triggered by a periodic coefficient.
Let $a\colon\Rd\times\Rd\to\mathbb R$ be $1$-periodic in both variables with $0<\alpha\le a\le\beta$, and $Q=(0,1)^d$. We consider the (squared) Gagliardo seminorms with oscillating coefficient
\begin{equation}\label{1}
\iint_{\Omega\times\Omega} a\Big(\frac x\eps,\frac y\eps\Big)\frac{|u(x)-u(y)|^2}{|x-y|^{d+2s}}\,dx\,dy ,
\end{equation}
with $\Omega$ either $\Rd$ or a Lipschitz open subset of $\Rd$. Two parameters compete: the period $\eps$ and the exponent $s=s(\eps)$. In order to make the results clearer we only treat the case of power $2$, and not of seminorms in $W^{s,p}$ with $p\neq 2$.

\smallskip
We can state this problem in the terminology of $\Gamma$-convergence.
For fixed $s\in(0,1)$, or for $s=s(\e)\to s_0\in(0,1)$ (in which case $s$ is replaced by $s_0$ in the limit), the $\Gamma$-limit with respect to the strong $L^2$ topology is
\begin{equation}
\overline a\iint_{\Omega\times\Omega}\frac{|u(x)-u(y)|^2}{|x-y|^{d+2s}}\,dx\,dy,\qquad\hbox{with}\qquad \overline a=\iint_{Q^2}a(x,y)\,dx\,dy .
\end{equation}
The reason is elementary: the contribution of pairs close to the diagonal is negligible, and, away from the diagonal, strong convergence together with a Riemann--Lebesgue argument replaces the coefficient by its average. In particular, there is no homogenization effect. Interesting behaviour can only appear when $s$ approaches an endpoint, either $1$ or $0$. In the homogeneous case, this asymptotic analysis was initiated in a seminal paper by Bourgain, Brezis, and Mironescu \cite{BBM}, who highlighted how, when $a$ is a constant, in the limit as $s\to 1$, the relevant interactions are those concentrating on the diagonal, for which the difference quotient in \eqref{1} can be approximated by a gradient. The opposite case, where $s\to 0$, has been analyzed in a paper by Maz'ya and Shaposhnikova \cite{MS}, who remarked that, in the case $\Omega=\Rd$, the relevant interactions are those with one of the points tending to infinity, producing a multiple of the squared $L^2$-norm. More precisely, if $a\equiv1$, for which the functional in \eqref{1} is the square of the  Gagliardo seminorm $[u]_{H^s}$,
two classical results dictate the normalization and the nature of the limit.

\smallskip\noindent {\bf The Maz'ya--Shaposhnikova formula} for $s\to0$. In this case, the meaningful set of integration is $\Rd$, for which we have $$\ s\,[u]^2_{H^s(\Rd)}\to\sig\|u\|^2_{L^2(\Rd)},$$ with $\sig$ being the $(d-1)$-dimensional surface area of the unit sphere in $\Rd$. The limit is a \emph{zero-order} functional that arises from interactions between points in the support of $u$ and ``points at infinity'', and can also be seen as a $\Gamma$-limit with respect to the weak $L^2$-convergence \cite{PoDeKuNiCr21}. For refinements and extensions of this formula, also to fractional perimeters and to metric measure spaces, we refer to \cite{DavoliDiFrattaGiorgioPinamonti25,DavFanPic26,DFPV,HanPinamonti21,HanPinamontiXuZambanini25}.

\smallskip\noindent{\bf The Bourgain--Brezis--Mironescu formula}  for $s\to1$. In this case, we have  $$\ (1-s)\,[u]^2_{H^s(\Omega)}\to\frac{\sig}{2d}\|\nabla u\|^2_{L^2(\Omega)}.$$ The limit is a \emph{first-order} functional, arising from interactions concentrating on the diagonal. In this case, the same result is a $\Gamma$-limit with respect to the strong convergence in $L^2(\Omega)$ \cite{ponce}.

\bigskip
In a recent pair of papers \cite{BMSP,BPSP} it has been shown that when oscillations are acting, contrary to the cases studied in \cite{BBM,MS}, we always have a contribution of interactions at a finite scale. Their effect is different in the case $s\to 0$, when they compete with oscillations at infinite scale, and in the case $s\to 1$, when their contribution is determined by the optimization of microscopic oscillations. In both cases, a parameter $\lambda$ appears, which describes the distribution of the energetic contribution between the various interactions. We provide a comparative description of the two results and propose a unified statement as an expansion by $\Gamma$-convergence.

\section{Determination of the relevant scales}
According to the two cases described above for $a\equiv1$, one scales the double integral in \eqref{1} by $s$ if $s\to 0$ and by $1-s$ as $s\to 1$. By integral-representation results, and since the functionals are quadratic forms, one expects the limit to be of the form
\begin{equation}\label{2}
b\int_{\Rd} |u|^2dx,\qquad\hbox{ and }\qquad \int_\Omega \langle B\nabla u,\nabla u\rangle dx
\end{equation}
in the respective cases.
The form of the coefficient $b$ and the matrix $B$ will be determined by the mutual behaviour of $\e$ and $s$. From now on, we choose $s=s(\e)$ and drop the dependence on $\e$.

\smallskip
The fundamental ansatz is that the behaviour of the functionals can be understood by looking at three types of contributions:

\smallskip
1) \emph{microscopic interactions}; that is, when $|x-y|\ll\eps$,

\smallskip
2) \emph{mesoscopic interactions}; that is, when $\e \lesssim|x-y|\lesssim 1$;

\smallskip
3) \emph{macroscopic interactions}; that is, when $|x-y|\gg 1$.

\bigskip\noindent
\emph{Scaling parameter as $s\to 0$}. The parameter $\e$ enters the kernel mass at the mesoscopic and macroscopic levels  as
\begin{equation}\label{4}
2s\int_{\Rd\setminus B_\e}\frac{1}{|\xi|^{d+2s}}d\xi=
2s\sig\int_\e^{+\infty}t^{-1-2s}\,dt=\sig \e^{-2s}
\end{equation}
(the factor $2$ is due to the symmetry of the interactions), which gives the relevant parameter 
$$
\lambda=\lim_{\eps\to0}\eps^{-2s}\in[1,+\infty].
$$
With this notation, the asymptotic mass at the mesoscopic level is 
\begin{equation}\label{4b}
\lim_{\e\to 0} 2s\int_{B_R\setminus B_\e}\frac{1}{|\xi|^{d+2s}}d\xi=
\lim_{\e\to 0}2s\sig\int_\e^{R}t^{-1-2s}\,dt=\sig (\lambda-1)
\end{equation}
for all $R>0$.

\medskip\noindent
\emph{Scaling parameter as $s\to 1$}.  The parameter $\e$ enters the kernel only at the microscopic level, where we can suppose $|u(x+\xi)-u(x)|\sim|\langle \nabla u(x),\xi\rangle|$.  The kernel contribution at $x$ is then of the order
\begin{equation}\label{5}
(1-s)\int_{B_\e}\frac{|\langle \nabla u(x),\xi\rangle|^2}{|\xi|^{d+2s}}d\xi=
(1-s)\frac{\sig}d|\nabla u(x)|^2\int_0^\e t^{1-2s}\,dt=\frac{\sig}{2d}|\nabla u(x)|^2\e^{2(1-s)},
\end{equation}
which gives the relevant parameter 
$$
\lambda=\lim_{\eps\to0}\eps^{2(1-s)}\in[0,1].
$$
With this notation, the asymptotic mass at the mesoscopic level is 
\begin{equation}\label{5b}
\lim_{\e\to 0} (1-s)\int_{B_R\setminus B_\e}\frac{|\langle \nabla u(x),\xi\rangle|^2}{|\xi|^{d+2s}}d\xi=\frac{\sig}{2d}|\nabla u(x)|^2(1-\lambda)
\end{equation}
for all $R>0$. Thanks to the integrability properties of the kernel at infinity, this argument may be applied at a scale $R=R_\e$ such that we may still assume 
$|u(x+\xi)-u(x)|\sim|\langle \nabla u(x),\xi\rangle|$, and at the same time interactions at a scale larger than $R_\e$ are negligible.


\section{The results}
In order to state the results, we need to define some quantities obtained from the coefficient $a$:
\smallskip

1) the {\em average} of $a$
$$
\overline a=\iint_{Q\times Q} a(x,y)dxdy.
$$
This quantity appears in {\em mesoscopic} interactions;

\smallskip
2) the {\em contracted} average of $a$ and its {\em harmonic mean} $\underline a$ given by
\begin{equation}\label{6}
\ac (x)=\int_Q (a(x,y)+a(y,x))dy, \qquad \frac1{\underline a}=\int_Q\frac1{\ac (x)}dx
\end{equation}
The coefficient $\ac (x)$ describes the averaged {\em macroscopic} interaction of $x$ with `points at infinity'. The harmonic mean is also characterized by the homogenization formula
$$
\Gamma\hbox{-}\lim_{\e\to 0} \int_\Omega\ac \Big(\frac x\e\Big)|u|^2dx=
\underline a\int_\Omega|u|^2dx,
$$
where the $\Gamma$-limit is taken with respect to the weak convergence in $L^2(\Omega)$;

\smallskip
3) the {\em diagonal term} $a(x,x)$ appears in the {\em microscopic interactions}. In order for this term to be well defined, for simplicity we additionally assume that $a$ is continuous. Correspondingly, we have a {\em homogenized matrix} $A_{\rm hom}$ characterized by 
$$
\Gamma\hbox{-}\lim_{\e\to 0} \int_\Omega a\Big(\frac x\e,\frac x\e\Big)|\nabla u|^2dx=
\int_\Omega \langle A_{\rm hom}\nabla u,\nabla u\rangle dx,
$$
where the $\Gamma$-limit is taken with respect to the weak and strong convergence in $L^2(\Omega)$.

\subsection{A Maz'ya--Shaposhnikova-type Homogenization Theorem}

The asymptotic characterization as $s\to 0$ is described by the following result.
\smallskip
\begin{theorem}[MS Homogenization Theorem \cite{BPSP}]\label{MS-th}  Let  $s\to0$, let $\eps^{-2s}\to\lambda\in[1,+\infty]$, and define
\begin{equation}\label{7}
a_\lambda=\inf_{\substack{\phi\in L^2(Q)\\ \int_Q\phi=0}}\Big\{\int_Q\ac(x)|1+\phi(x)|^2dx+(\lambda-1)\iint_{Q\times Q}a(x,y)|\phi(x)-\phi(y)|^2dx\,dy\Big\}
\end{equation} 
if $\lambda<+\infty$, and
\begin{equation}\label{8}
a_{+\infty}=\int_Q\ac(x)dx=2\overline a
\end{equation} 
if $\lambda=+\infty$.
Then we have 
\begin{equation}
\Gamma\hbox{-}\lim_{\e\to 0}s\iint_{\Rd\times\Rd}a\Big(\frac x\e,\frac y\e\Big)\frac{|u(x)-u(y)|^2}{|x-y|^{d+2s}}dxdy =\frac{\sig}{2}\,a_\lambda\,\|u\|_{L^2}^2=: F^{\rm MS}_\lambda(u)
\end{equation}
 in the sequential \emph{weak} $L^2(\Rd)$ topology.
 \end{theorem}

 \smallskip 
 The optimal $\phi$ in \eqref{7} gives the so-called  {\em corrector}: extending it by periodicity, given a Lipschitz bounded open set $A\subset\Rd$, the functions $u_\e(x)=\chi_A(x) (1+\phi(\frac x\e))$ are such that if $F_\e$ denotes $s$ times the functional in \eqref{1} with $\Omega=\Rd$, $F_\e(u_\e)\to F^{\rm MS}_\lambda(\chi_A)$.
 This result shows that in this regime, oscillations are felt both in the mesoscopic and in the macroscopic interactions, while the effect of microscopic interactions is negligible. The optimization of oscillations at both scales gives $a_\lambda$.   The prefactor $\frac{\sig}{2}$ in the $\Gamma$-limit and the factor $\lambda -1$ in the definition of $a_\lambda$ are heuristically explained by \eqref{4} and \eqref{4b}. 
 Note that for $a\equiv1$, we have $a_\lambda=2$, and one recovers $\sig\|u\|^2_{L^2}$.

\smallskip
The map $\lambda\mapsto a_\lambda$ is nondecreasing, continuous, and concave in $[1,+\infty]$. 
The case $\lambda=+\infty$ corresponds to \eqref{8}, which can be read as \eqref{7} with the factor $\lambda-1$ understood as forcing the corrector to vanish. 
For $\lambda=1$, $a_1=\underline a$ is the harmonic mean of $\ac $  as defined in \eqref{6}. In this regime, the mesoscopic interactions are negligible, and the functionals are 
equivalent to
$$
F^1_\e(u)=\frac{\sig}{2}\int_{\Rd} \ac \Big(\frac x\e\Big)|u|^2dx;
$$
that is, we can first average `interactions at infinity' and then homogenize with respect to the weak $L^2$ topology. 

\subsection{A Bourgain--Brezis--Mironescu-type Homogenization Theorem}

The asymptotic characterization as $s\to 1$ is the following.
\smallskip

\begin{theorem}[BBM Homogenization Theorem \cite{BMSP}]\label{BBM-th} Let $\Omega$ be a Lipschitz open set, and let $a$ be continuous. Let $s\to1$, let $\eps^{2(1-s)}\to\lambda\in[0,1]$, and define
\begin{equation}\label{10}
A_\lambda= (1-\lambda)\overline a\,{\rm Id}+\lambda A_{\rm hom},
\end{equation} 
where {\rm Id} denotes the identity matrix.
Then we have 
\begin{equation}\label{11}
\Gamma\hbox{-}\lim_{\e\to 0}(1-s)\iint_{\Omega\times\Omega}a\Big(\frac x\e,\frac y\e\Big)\frac{|u(x)-u(y)|^2}{|x-y|^{d+2s}}dxdy =\frac{\sig}{2d}\int_\Omega \langle A_{\lambda}\nabla u,\nabla u\rangle dx=: F^{\rm BBM}_\lambda(u)
\end{equation}
equivalently in the weak or strong $L^2(\Omega)$ topologies.
\end{theorem}
 
 \smallskip
Formula \eqref{10} highlights a different behaviour than in the case $s\to 0$. Again, interactions are spread across two levels: the microscopic and the mesoscopic. Contrary to the case $s\to 0$, the homogenization of the oscillations now only takes place at the microscopic level, while their average contributes at the mesoscopic level. As a result, the limit is an exact convex combination of the averaged energies (that is, the limit at $s>0$ fixed) and the homogenized energies in $H^1(\Omega)$ with coefficient $a(x,x)$. The correctors in this case are the same as for the homogenization problem. More precisely, with fixed $\xi\in\mathbb R^d$, one takes $\phi_\xi$ $1$-periodic and minimizing
$$
\langle A_{\rm hom}\xi,\xi\rangle= \min\Big\{ \int_Q a(y,y)|\xi+\nabla \phi|^2dy: \phi\in H^1_{\rm per}(Q)\Big\};
$$
then, if we set $u^\xi(x)=\langle\xi,x\rangle$, 
the functions $\langle\xi,x\rangle+ \e\phi_\xi(\frac x\e)$ are such that
$$
\lim_{\e\to 0} (1-s)\iint_{\Omega\times\Omega}a\Big(\frac x\e,\frac y\e\Big)\frac{|u^\xi(x)+\e\phi_\xi(\frac x\e)-u^\xi(y)-\e\phi_\xi(\frac y\e)|^2}{|x-y|^{d+2s}}dxdy = F^{\rm BBM}_\lambda(u^\xi).
$$
Formula \eqref{10} shows that the form of $\phi_\xi$ does not influence the asymptotic behaviour of the mesoscopic part of the energies. The convex combination is heuristically explained by \eqref{5b}, where the relation
$|u(x+\xi)-u(x)|\sim|\langle \nabla u(x),\xi\rangle|$ is justified by a Riemann--Lebesgue type lemma.
Note that for $a\equiv1$, \eqref{11} is the BBM formula. 

\bigskip\noindent
{\bf Remark}~({\em separation-of-scale regimes}).
In both cases, $s\to 0^+$ and $s\to 1^-$, the extreme values of the corresponding $\lambda$ determine the regimes where we have a \emph{separation of scales}; that is, the corresponding limits can be (at least, formally) interpreted as iterated $\Gamma$-limits.

\smallskip
{\em Case $s\to 0$ and  $|\log\e|\ll\frac1s$}. In this case, we have $\lambda=1$ and $a_{1}=\underline a$. This also corresponds to taking the limit first as $s\to0$ and then as $\e\to 0$. Indeed, formally,  with fixed $\e$, the first gives
$$
\Gamma\hbox{-}\lim_{s\to 0}s\iint_{\Rd\times\Rd} a\Big(\frac x\e,\frac y\e\Big)\frac{|u(x)-u(y)|^2}{|x-y|^{d+2s}}dxdy =\frac{\sig}{2}\int_{\Rd} \ac\Big(\frac x\e\Big)|u(x)|^2dx,
$$
and then we take the $\Gamma$-limit as $\e\to 0$. This gives as coefficient the harmonic mean of $\ac$; that is, $\underline a$.

\smallskip
{\em Case $s\to 0$ and $|\log\e|\gg\frac1s$}. In this case,  $\lambda=+\infty$, and we have $a_{+\infty}=2\overline a$, which also corresponds to 
taking the limit first as $\e\to0$ and then as $s\to 0$. Indeed, we have seen that, with fixed $s>0$
$$
\Gamma\hbox{-}\lim_{\e\to 0}s\iint_{\Rd\times\Rd} a\Big(\frac x\e,\frac y\e\Big)\frac{|u(x)-u(y)|^2}{|x-y|^{d+2s}}dxdy =s\,\overline a\iint_{\Rd\times\Rd} \frac{|u(x)-u(y)|^2}{|x-y|^{d+2s}}dxdy. 
$$ 
We can then take the limit as $s\to 0$ and apply the Maz'ya--Shaposhnikova result to obtain $2\overline a$.

\smallskip
{\em Case $s\to 1$ and $|\log\e|\gg\frac1{1-s}$}. In this case,  $\lambda=0$, and we have $A_0=\overline a\,{\rm Id}$, which also corresponds to 
taking the limit first as $\e\to0$ and then as $s\to 1$. Indeed, by the observation of the previous case, this is equivalent to taking
$$
\Gamma\hbox{-}\lim_{s\to 1}\overline a(1-s)\iint_{\Omega\times\Omega} \frac{|u(x)-u(y)|^2}{|x-y|^{d+2s}}dxdy, 
$$ 
which gives $\overline a\frac{\sig}{2d}\|\nabla u\|^2_{L^2(\Omega)}$ by the BBM formula, corresponding to $\overline a\,{\rm Id}$.

\smallskip
{\em Case $s\to 1$ and $|\log\e|\ll\frac1{1-s}$}. In this case,  $\lambda=1$, and we have $A_1=A_{\rm hom}$. Formally, this corresponds to first  applying the BBM formula at fixed $\e>0$ and 
obtaining
$$
\Gamma\hbox{-}\lim_{s\to 1}(1-s)\iint_{\Omega\times\Omega} a\Big(\frac x\e,\frac y\e\Big)\frac{|u(x)-u(y)|^2}{|x-y|^{d+2s}}dxdy =\frac{\sig}{2d}\int_{\Omega} a\Big(\frac x\e,\frac x\e\Big)|\nabla u|^2dx,
$$
and then homogenizing the latter energy. Note that the concentration on the diagonal gives the coefficient with $x=y$ in the right-hand side expression.

\section{Expansion by $\Gamma$-convergence}
The above results can be summarized in terms of expansion by $\Gamma$-convergence \cite{BT}. 
In our context, the family $G_{\e,s}$ is a {\em $\Gamma$-expansion} (at order $0$) of the family $F_{\e,s}$ 
for $\e\to 0$ and $0<s<1$ simply if we have
$$
\Gamma\hbox{-}\lim_{h\to+\infty} F_{\e_h,s_h}= \Gamma\hbox{-}\lim_{h\to+\infty} G_{\e_h,s_h}
$$
for every pair of sequences $\e_h\to 0$ and $s_h\to s_0\in[0,1]$ such that both $\Gamma$-limits exist.
This equivalence is designed so that minimum problems involving $F_{\e,s}$ are uniformly approximated by the
corresponding minimum problems with $G_{\e,s}$, under an equi-coerciveness assumption. 
Of course, this fact is meaningful if $G_{\e,s}$ are functionals that are simpler to handle.
In our case, the simplification is that $G_{\e,s}$ are homogeneous functionals; that is, their integrands are independent of $x$. Note that, 
since the definition is given in terms of $\Gamma$-limits, it is natural that the domain of $F_{\e,s}$ may differ from that of $G_{\e,s}$. In this particular case, we can think of all functionals as defined in $L^2(\Rd)$. We consider the weak convergence in $L^2(\Rd)$, but we can equivalently consider the strong convergence in $L^2(\Rd)$ if $s\ge s_0>0$.

For the convenience of exposition, in the following statement, the functionals $G_{\e,s}$ are defined by distinguishing three regimes, which are to be understood as conditions on sequences $(\e_h,s_h)$ with $\e_h\to0$; that is, $G_{\e_h,s_h}$ is given by \eqref{13}, \eqref{14}, or \eqref{15} according to the asymptotic behavior of the sequence. In the remark following the statement, we specify how $G_{\e,s}$ can be defined for each pair $(\e,s)$, as required by the definition of $\Gamma$-expansion.

\begin{theorem}[$\Gamma$-expansion of oscillating Gagliardo energies]\label{thm3}
Let $a$, $a_\lambda$, and $A_\lambda$ be defined as in the previous part of the paper. For all $s\in(0,1)$ and $\e>0$ define
\begin{equation}\label{12}
F_{\e,s}(u)=s(1-s)\iint_{\Rd\times\Rd} a\Big(\frac x\eps,\frac y\eps\Big)\frac{|u(x)-u(y)|^2}{|x-y|^{d+2s}}\,dx\,dy ,
\end{equation}
for $u\in H^s(\Rd)$, and $G_{\e,s}(u)$ as follows (for notational ease, we set ${\bf a}(\lambda)=a_\lambda$ and ${\bf A}(\lambda)=A_\lambda$):

\smallskip
{\rm 1)} if $|\log\e|\gg \frac1{s(1-s)}$ 
\begin{equation}\label{13}
G_{\e,s}(u)=s(1-s)\overline a\iint_{\Rd\times\Rd} \frac{|u(x)-u(y)|^2}{|x-y|^{d+2s}}\,dx\,dy
\end{equation}
for $u\in H^s(\Rd)$;

\smallskip
{\rm 2)} if $|\log\e|= O(\frac1s)$
\begin{equation}\label{14}
G_{\e,s}(u)={\bf a}(\eps^{-2s})\frac{\sig}{2}\int_{\Rd}|u|^2dx
\end{equation}
if $u\in L^2(\Rd)$;

\smallskip
{\rm 3)} if $|\log\e|= O(\frac1{1-s})$
\begin{equation}\label{15}
G_{\e,s}(u)=\frac{\sig}{2d}\int_{\Rd} \langle {\bf A}(\e^{2(1-s)})\nabla u,\nabla u\rangle dx
\end{equation}
if $u\in H^1(\Rd)$.

\smallskip

Then $G_{\e,s}$ is a $\Gamma$-expansion (at order $0$) of the family $F_{\e,s}$ 
for $\e\to 0$ and $0<s<1$ with respect to the weak convergence in $L^2(\Rd)$. 
\end{theorem}

\smallskip\noindent
{\bf Remark}~({\em definition of $G_{\e,s}$ for each pair $(\e,s)$}). Conditions (1)--(3) are meaningful only along sequences, and do not cover sequences along which, for instance, $s|\log\e|$ oscillates between bounded and unbounded values. In order to have $G_{\e,s}$ defined for each pair $(\e,s)$, it suffices to fix $\eta_\e\to0$ with $\eta_\e|\log\e|\to+\infty$ (e.g., $\eta_\e=|\log\e|^{-1/2}$), and to define $G_{\e,s}$ by \eqref{14} if $s\le\eta_\e$, by \eqref{15} if $s\ge 1-\eta_\e$, and by \eqref{13} otherwise. Indeed, all sequences with $|\log\e|=O(\frac1s)$ eventually satisfy $s\le\eta_\e$; moreover, if $s\le\eta_\e$ then $s\to0$, and the functionals in \eqref{14} have the correct $\Gamma$-limit also when $\e^{-2s}\to+\infty$, by the continuity of ${\bf a}$ at $+\infty$. Conversely, if $s>\eta_\e$ and $s\to 0$ then $\e^{-2s}\to+\infty$, a case in which the functionals in \eqref{13} and \eqref{14} have the same $\Gamma$-limit $\sig\,\overline a\,\|u\|^2_{L^2}$. A symmetric argument holds for $s$ close to $1$, using the continuity of ${\bf A}$ at $0$.

\bigskip
\noindent{\em Proof of Theorem {\rm\ref{thm3}}.} If $|\log\e|\gg ({s(1-s)})^{-1}$, the claim follows from the convergence for $s$ to $s_0\in(0,1)$, from Theorem \ref{MS-th} in the case $\lambda=+\infty$, and from Theorem \ref{BBM-th} (with $\Omega=\Rd$) in the case $\lambda=0$, noting that the last two results, applied with the constant coefficient $\overline a$ in the place of $a$, give the same $\Gamma$-limits for the functionals in \eqref{13}. Theorem \ref{MS-th} in the case $\lambda<+\infty$ proves the claim in the second case, and Theorem \ref{BBM-th} in the case $\lambda>0$  (and with $\Omega=\Rd$) proves the claim in the third case. In the last two cases, we use the continuous dependence on $\lambda$ of the coefficient ${\bf a}(\lambda)$ and of the matrix ${\bf A}(\lambda)$. With $G_{\e,s}$ defined as in the remark above, the claim follows from the same arguments, upon extraction of subsequences.\hfill $\square$

\bigskip\noindent
{\bf Remark}~({\em extension of the expansion to $0$ and $1$}). The $\Gamma$-expansion result can be extended to $s\in[0,1]$ by defining
$$
F_{\e,0}(u)=\frac{\sig}{2}\int_{\Rd} \ac \Big(\frac x\eps\Big)|u(x)|^2\,dx ,
\qquad
F_{\e,1}(u)=\frac{\sig}{2d}\int_{\Rd}a\Big(\frac x\eps,\frac x\eps\Big) |\nabla u|^2\,dx ,
$$
with the corresponding $G_{\e,0}$ and $G_{\e,1}$, independent of $\e$, 
defined by
$$
G_{\e,0}(u)=\frac{\sig}{2}\underline a\int_{\Rd}|u(x)|^2\,dx ,
\qquad
G_{\e,1}(u)=\frac{\sig}{2d}\int_{\Rd}\langle A_{\rm hom}\nabla u,\nabla u\rangle\,dx;
$$
that is, $G_{\e,0}$ is \eqref{14} evaluated at $\lambda=1$ (or, equivalently, when $\e^{-2s}\to1$), and $G_{\e,1}$ is \eqref{15} at $\lambda=1$ (or, equivalently, when $\e^{2(1-s)}\to 1$).
Indeed, the only cases left out from the two homogenization theorems are when $s=0$ or $s=1$ and $\e\to 0$,
which are ensured by the classical homogenization theorems already mentioned.

\bigskip\noindent
{\bf Acknowledgements.} 
 Andrea Braides acknowledges funding by the Excellence Project MatMod@TOV (Grant Agreement No.~CUP\_E83C23000330006) of the
Department of Mathematics of the University of Rome Tor Vergata. The authors are members of GNAMPA, INdAM.

\bigskip
\noindent\textbf{Authors' addresses}

\smallskip
\noindent Andrea Braides, Dipartimento di Matematica, Universit\`a di Roma Tor Vergata,
Via della Ricerca Scientifica 1, 00133 Roma, Italy. E-mail: \texttt{[braides@mat.uniroma2.it]}

\smallskip
\noindent S. Mannella, M. Picerni, A.J.A. Schiavoni Piazza, SISSA, Via Bonomea 265,
34136 Trieste, Italy. E-mail: \texttt{[smannell@sissa.it]}, \texttt{[mpicerni@sissa.it]}, \texttt{[aschiavo@sissa.it]}

\end{document}